\documentclass[12pt]{article}

      \usepackage{latexsym}
         \usepackage[reqno, namelimits, sumlimits]{amsmath}
         \usepackage{amssymb, amsfonts}
         \usepackage{amsthm}

 \newtheorem{theorem}{Theorem}[section]
 
 \newtheorem{definition}[theorem]{Definition}

 \newtheorem{remark}[theorem]{Remark}

 \newtheorem{pro}[theorem]{Proposition}

\title{On Type I Blowups of Suitable Weak Solutions to Navier-Stokes Equations near Boundary.
}

 \author{ G.~Seregin\footnote{Oxford University, UK, and St Petersburg Department of Steklov Mathematical Institute, RAS, Russia, email: seregin@maths.ox.ac.uk}
}

\date{}
\begin{document}
\maketitle
\begin{abstract} In this note, boundary Type I blowups of suitable weak solutions to the Navier-Stokes equations are discussed. In particular, it has been shown that, under certain assumptions, the existence of non-trivial mild bounded ancient solutions in half space leads to the existence of suitable weak solutions with Type I blowup on the boundary.
 
 \end{abstract}

  \section{Introduction}
 \setcounter{equation}{0}

 The aim of the note is to study conditions under which solutions to the Navier-Stokes equations undergo Type I blowups on the boundary.
 
 Consider the classical Navier-Stokes equations  \begin{equation}
 \label{NSE}
 \partial_tv+v\cdot\nabla v-\Delta v=-\nabla q,\qquad {\rm div}\,v=0
 	\end{equation}
 	in the space time domain $Q^+=B^+\times ]-1,0[$, where $B^+=B^+(1)$ and $B^+(r)=\{x\in \mathbb R^3: |x|<r,\,x_3>0\}$ is a half ball of radius $r$ centred at the origin $x=0$. It is supposed that $v$ satisfies the homogeneous Dirichlet boundary condition
\begin{equation}\label{bc}
 	v(x',0,t)=0	
 	\end{equation}
for all $|x'|<1$ and $-1<t<0$. Here, $x'=(x_1,x_2)$ so that $x=(x',x_3)$ and $z=(x,t)$.

Our goal is to understand  how to determine whether or not the origin $z=0$ is a singular point of the velocity field $v$. We say that $z=0$ is a regular point of $v$ if there exists $r>0$ such that  $v\in L_\infty(Q^+(r))$ where $Q^+(r)=B^+(r)\times ]-r^2,0[$. It is known, see \cite{S3}   and \cite{S2009},
 that the velocity $v$ is H\"older continuous in a parabolic vicinity of $z=0$ if $z=0$ is a regular point. However, further smoothness even in spatial variables does not follow in the regular case, see \cite{Kang2005} and \cite{SerSve10} for counter-examples.

  The class of solutions to be  studied  is as follows.
  
  \begin{definition}\label{sws} A pair of functions $v$ and $q$ is called a suitable weak solution to the Navier-Stokes equations in $Q^+$ near the boundary if and only if the following conditions hold:
 \begin{equation}\label{class}
 v\in L_{2,	\infty}(Q^+)\cap L_2(-1,0;W^1_2(Q^+)),\qquad q\in L_\frac 32(Q^+);
 \end{equation}
 $v$ and $q$ satisfy equations (\ref{NSE}) and boundary condition (\ref{bc});
 	$$\int\limits_{B^+}\varphi(x,t)|v(x,t)|^2dx+2\int\limits_{-1}^t\int\limits_{B^+}\varphi|\nabla v|^2dxdt\leq 	$$
 	\begin{equation}\label{energy_inequality}
 	\leq	\int\limits_{-1}^t\int\limits_{B^+}(|v|^2(\partial_t\varphi+\Delta\varphi)+v\cdot\nabla v(|v|^2+2q))dxdt 	\end{equation}
 for all non-negative functions $\varphi\in C^\infty_0(B\times]-1,1[)$ such that $\varphi|_{x_3=0}=0$.	
 \end{definition}
 
 In what follows, some statements will be expressed in   terms of scale invariant quantities (invariant with respect to the Navier-Stokes scaling: $\lambda v(\lambda x,\lambda^2 t)$ and $\lambda ^2q(\lambda x,\lambda^2 t)$). Here, they are:
 $$A(v,r)=\sup\limits_{-r^2<t<0}\frac 1r\int\limits_{B^+(r)}|v(x,t)|^2dx, \qquad
 E(v,r)=\frac 1r\int\limits_{Q^+(r)}|\nabla v|^2dz,$$$$
 C(v,r)=\frac 1{r^2}\int\limits_{Q^+(r)}|v|^3dz,\qquad 
 D_0(q,r)=\frac 1{r^2}\int\limits_{Q^+(r)}|q-[q]_{B^+(r)}|^\frac 32dz,  $$
 $$D_2(q,r)=\frac 1{r^\frac {13}8}\int\limits^0_{-r^2}\Big(\int\limits_{B^+(r)}|\nabla q|^\frac {12}{11}dx\Big)^\frac {11}8dt,$$
 where 
 $$[f]_\Omega=\frac 1{|\Omega|}\int\limits_\Omega fdx.$$

We also introduce the following values:
$$g(v):=\min\{\sup\limits_{0<R<1}A(v,R), \sup\limits_{0<R<1}C(v,R),\sup\limits_{0<R<1}E(v,R)\}
$$
and, given $r>0$,
$$G_r(v,q):=$$$$\max\{\sup\limits_{0<R<r}A(v,R), \sup\limits_{0<R<r}C(v,R),\sup\limits_{0<R<r}E(v,R),\sup\limits_{0<R<r}D_0(q,R)\}.
$$

Relationships between $g(v)$ and $G_1(v,q)$ is described in the following proposition.

\begin{pro}\label{boundednesstheorem}
Let $v$ and $q$ be a suitable weak solution to the Navier-Stokes equations in $Q^+$ near the boundary. Then, $G_1$ is bounded if and only if $g$ is bounded.	
\end{pro}

 If $z=0$ is a singular point of $v$ and $g(v)<\infty$, then $z=0$ is called a Type I singular point or a Type I blowup point.

Now, we are ready to state the main results of the paper.

\begin{definition}
	\label{leas}
A  function $u:Q^+_-:=\mathbb R^3_+\times]-\infty,0[\,\to\mathbb R^3$ is called a local energy ancient solution if there exists a function $p:Q_-^+\to\mathbb R$ such that  the pair $u$ and $p$ is a suitable weak solution in $Q^+(R)$ for any $R>0$. Here, $\mathbb R^3_+:=\{x\in \mathbb R^3:\,x_3>0\}$.
\end{definition} 
 
\begin{theorem}\label{local energy ancient solution} There exists a suitable weak solution $v$ and $q$ with Type I blowup at the origin $z=0$ if and only if there exists a non-trivial local energy ancient solution $u$ such that $u$ and the corresponding pressure $p$ have the following prosperities:
\begin{equation}\label{scalequatities}
G_\infty(u,p):=\max\{\sup\limits_{0<R<\infty} A(u,R), \sup\limits_{0<R<\infty}E(u,R),$$$$\sup\limits_{0<R<\infty}C(u,R),\sup\limits_{0<R<\infty}D_0(p,R)\}<\infty	
\end{equation}
and
\begin{equation}\label{singularity}
\inf\limits_{0<a<1}C(u,a)\geq \varepsilon_1>0.	
\end{equation}
 	
 \end{theorem}
\begin{remark}\label{singType1}
According to (\ref{scalequatities}) and (\ref{singularity}), the origin $z=0$ is Type I blowup of the velocity $u$.	
\end{remark}

 There is another way to construct a suitable weak solution with Type I blowup. It is motivated by the recent result in \cite{AlBa18} for the interior case. Now, the main object is related to the so-called mild bounded ancient solutions in a half space, for details see \cite{SerSve15} and \cite{BaSe15}.
\begin{definition}\label{mbas}
A bounded	function $u$ is a mild bounded ancient solution if and only if there exists a pressure $p=p^1+p^2$, where the even extension of $p^1$ in $x_3$ to the whole space $\mathbb R^3$ 
is a $L_\infty(-\infty,0;BMO(\mathbb R^3))$-function,
$$\Delta p^1={\rm divdiv}\,u\otimes u$$
in $Q^+_-$ with $p^1_{,3}(x',0,t)=0$, and $p^2(\cdot,t)$ is a harmonic function in $\mathbb R^3_+$, whose gradient satisfies the estimate
$$|\nabla p^2(x,t)|\leq \ln (2+1/x_3)$$
for all $(x,t)\in Q^+_-$ and has the property
$$\sup\limits_{x'\in \mathbb R^2}|\nabla p^2(x,t)|\to 0
$$ as $x_3\to \infty$; functions $u$ and $p$ satisfy:
$$\int\limits_{Q^+_-}u\cdot \nabla qdz=0$$
for all $q\in C^\infty_0(Q_-:=\mathbb R^3\times ]-\infty,0[)$ and, for any $t<0$,
$$\int\limits_{Q^+_-}\Big(u\cdot(\partial_t\varphi+\Delta\varphi)+u\otimes u:\nabla \varphi+p{\rm div}\,\varphi\Big)dz=0
$$ for and $\varphi\in C^\infty_0(Q_-)$ with $\varphi(x',0,t)=$ for all $x'\in \mathbb R^2$.
\end{definition}

As it has been shown in \cite{BaSe15}, any mild bounded ancient solution $u$ in a half space is infinitely smooth up to the boundary and $u|_{x_3}=0$.
\begin{theorem}\label{mbas_type1}
Let $u$ be a mild bounded ancient solution such that $|u|\leq 1$ and $|u(0,a,0)|=1$ for a positive number $a$ and such that
 (\ref{scalequatities}) holds. Then there exists a suitable weak solution in $Q^+$ having Type I blowup at the origin $z=0$.
\end{theorem}
{\bf Acknowledgement} The work is supported by the grant RFBR 17-01-00099-a.

 \section{Basic Estimates}
 \setcounter{equation}{0}
 
 In this section, we are going to state and prove certain basic estimates for arbitrary suitable weak solutions near the boundary.

For our purposes, the main estimate of the convective term can be derived  as follows. First, we apply H\"older inequality in spatial variables:
$$\||v||\nabla v|\|_{\frac {12}{11},\frac 32,Q_+(r)}^\frac 32=\int\limits^0_{-r^2}\Big(\int\limits_{B_+(r)}
|v|^\frac {12}{11}|\nabla v|^\frac {12}{11}dx\Big)^\frac{11}8dt\leq$$
$$\leq\int\limits^0_{-r^2}\Big(\int\limits_{B_+(r)}|\nabla v|^2dx\Big)^\frac 34\Big(\int\limits_{B_+(r)}|v|^\frac {12}5dx\Big)^\frac 58dt.
$$
Then, byy interpolation, since $\frac {12}5=2\cdot\frac 35+3\cdot\frac 25$, we find
$$\Big(\int\limits_{B_+(r)}|v|^\frac {12}5dx\Big)^\frac 58\leq \Big(\int\limits_{B_+(r)}|v|^2dx\Big)^\frac 38\Big(\int\limits_{B_+(r)}|v|^3dx\Big)^\frac 14. $$
So,
$$\||v||\nabla v|\|_{\frac {12}{11},\frac 32,Q_+(r)}^\frac 32\leq $$$$\leq \int\limits^0_{-r^2}\Big(\int\limits_{B_+(r)}|\nabla v|^2dx\Big)^\frac 34\Big(\int\limits_{B_+(r)}|v|^2dx\Big)^\frac 38\Big(\int\limits_{B_+(r)}|v|^3dx\Big)^\frac 14dt\leq $$
\begin{equation}\label{mainest}\leq \sup\limits_{-r^2<t<0}\Big(\int\limits_{B_+(r)}|v|^2dx\Big)^\frac 38\Big(\int\limits_{Q_+(r)}|\nabla v|^2dxdt\Big)^\frac 34\Big(\int\limits_{Q_+(r)}|v|^3dxdt\Big)^\frac 14\leq \end{equation}
$$\leq r^\frac 38r^\frac 34r^\frac 12A^\frac 38(v,r)E^\frac 34(v,r)C^\frac 14(v,r)
$$$$=r^\frac {13}8A^\frac 38(v,r)E^\frac 34(v,r)C^\frac 14(v,r).
$$

Two other estimates are well known and  valid for any $0<r\leq 1$:
\begin{equation}\label{multiple}
C(v,r)\leq cA^\frac 34(v,r)E^\frac 34(v,r)
	\end{equation} 
	 and 
\begin{equation}\label{embedding}
	D_0(q,r)\leq cD_2(q,r).
\end{equation}

Next, one more estimate immediately follows from the energy inequality (\ref{energy}) for a suitable choice of cut-off function $\varphi$:
\begin{equation}\label{energy}
A(v,\tau R)+E(v,\tau R)\leq c(\tau)\Big[C^\frac 23(v,R)+C^\frac 13(v,R)D_0^\frac 23(q,R)+C(v,R)\Big]	
\end{equation}
for any $0<\tau<1$ and for all $0<R\leq 1$.
 
 The last two estimates are coming out from the linear theory. Here, they are:
\begin{equation}\label{pressure}
D_2(q,r)\leq c \Big(\frac r\varrho\Big)^2\Big[D_2(q,\varrho)+E^\frac 34(v,\varrho)\Big]+$$$$+c\Big(\frac \varrho r\Big)^\frac {13}8A^\frac 38(v,\varrho)E^\frac 34(v,\varrho)C^\frac 14(v,\varrho)	
\end{equation}
for any $0<r<\varrho\leq 1$ and 
\begin{equation}\label{highder}
\|\partial_tv\|_{\frac {12}{11},\frac 32,Q^+(\tau R)}+\|\nabla^2v\|_{\frac {12}{11},\frac 32,Q^+(\tau R)}+\|\nabla q\|_{\frac {12}{11},\frac 32,Q^+(\tau R)}	\leq 
\end{equation}
$$\leq c(\tau)R^\frac {13}{12}\Big[D_0^\frac 23(q,R)+C^\frac 13(v,R)+E^\frac 12(v,R)+$$$$+(A^\frac 38(v,R)E^\frac 34(v,R)C^\frac 14(v,R))^\frac 23\Big]
$$
for any $0<\tau<1$ and for all $0<R\leq 1$.



Estimate (\ref{highder}) follows from bound (\ref{mainest}), from the local regularity theory for the Stokes equations (linear theory), see paper
\cite{S2009}, and from scaling. Estimate (\ref{pressure}) will be proven in the next section.

\section{Proof of (\ref{pressure})}
 \setcounter{equation}{0} 
 Here, we follows paper \cite{S3}. We let 
 $\tilde f=-v\cdot\nabla v$ and observe
 that
 \begin{equation}\label{weakerterm}
 \frac 1r\|\nabla v\|_{\frac {12}{11},\frac 32,Q^+(r)}\leq r^\frac {13}{12}E^\frac 12(v,r)
	\end{equation}
and, see  (\ref{mainest}),
\begin{equation}\label{righthand}
\|\tilde f\|_{\frac {12}{11},\frac 32,Q^+(r)}	\leq cr^\frac {13}{12}(A^\frac 38(v,r)E^\frac 34(v,r)C^\frac 14(v,r))^\frac 23.
\end{equation}
Next, we select a  convex domain with smooth boundary so that
$$B^+(1/2)\subset \tilde B\subset B^+$$ and, for $0<\varrho<1$, we let 
$$\tilde B(\varrho)=\{x\in \mathbb R^3: x/\varrho\in \tilde B\},\qquad \tilde Q(\varrho)=\tilde B(\varrho)\times ]-\varrho^2,0[.
$$
Now,  consider the following initial boundary value problem:
\begin{equation}\label{v1eq}
\partial_tv^1-\Delta v^1+\nabla q^1=\tilde f, \qquad {\rm div}\,v^1=0	
\end{equation} in $\tilde Q(\varrho)$ and 
\begin{equation}\label{v1ibc}
v^1=0	
\end{equation}
on parabolic boundary $\partial'\tilde Q(\varrho)$ of $\tilde Q(\varrho)$. It is also supposed that
$[q^1]_{\tilde B(\varrho)}(t)=0$ for all $-\varrho^2<t<0$.

Due to estimate ({\ref{righthand}) and due to the Navier-Stokes  scaling, 
a unique solution to problem (\ref{v1eq}) and (\ref{v1ibc}) satisfies the estimate
$$\frac 1{\varrho^2}\|v^1\|_{\frac {12}{11},\frac 32,\tilde Q(\varrho)}+\frac 1{\varrho}\|\nabla v^1\|_{\frac {12}{11},\frac 32,\tilde Q(\varrho)}
+\|\nabla^2 v^1\|_{\frac {12}{11},\frac 32,\tilde Q(\varrho)}+$$
\begin{equation}\label{v1est}
+	\frac 1{\varrho}\| q^1\|_{\frac {12}{11},\frac 32,\tilde Q(\varrho)}
+\|\nabla q^1\|_{\frac {12}{11},\frac 32,\tilde Q(\varrho)}\leq 
\end{equation}
$$\leq c\|\tilde f\|_{\frac {12}{11},\frac 32,\tilde Q(\varrho)}\leq c\varrho^\frac {12}{11}(A^\frac 38(v,\varrho)E^\frac 34(v,\varrho)C^\frac 14(v,\varrho))^\frac 23,$$
where a generic constant c is independent of $\varrho$.

Regarding $v^2=v-v^1$ and $q^2=q-[q]_{B_+(\varrho/2)}-q^1$, one can notice the following:\begin{equation}\label{v2eq}
\partial_tv^2-\Delta v^2+\nabla q^2=0, \qquad {\rm div}\,v^2=0		
\end{equation} in $\tilde Q(\varrho)$ and 
\begin{equation}\label{v2ibc}
v^2|_{x_3=0}=0.	
\end{equation}
As it was indicated in \cite{S2009}, functions $v^2$ and $q^2$ obey the estimate
\begin{equation}\label{v2q2est}\|\nabla^2 v^2\|_{9,\frac 32, Q^+(\varrho/4)}+\|\nabla q^2\|_{9,\frac 32, Q^+(\varrho/4)}\leq \frac c{\varrho^\frac {29}{12}}L,\end{equation}
where
$$L:=\frac 1{\varrho^2}\| v^2\|_{\frac {12}{11},\frac 32, Q^+(\varrho/2)}+\frac 1{\varrho}\| \nabla v^2\|_{\frac {12}{11},\frac 32, Q^+(\varrho/2)}+\frac 1{\varrho}\| q^2\|_{\frac {12}{11},\frac 32, Q^+(\varrho/2)}.$$
As to an evaluation of $L$, we have
$$L\leq
\Big[\frac 1{\varrho^2}\| v\|_{\frac {12}{11},\frac 32, Q^+(\varrho/2)}+\frac 1{\varrho}\| \nabla v\|_{\frac {12}{11},\frac 32, Q^+(\varrho/2)}+\frac 1{\varrho}\| q-[q]_{B^+(\varrho/)}\|_{\frac {12}{11},\frac 32, Q^+(\varrho/2)}+$$
$$+\frac 1{\varrho^2}\| v^1\|_{\frac {12}{11},\frac 32, Q^+(\varrho/2)}+\frac 1{\varrho}\| \nabla v^1\|_{\frac {12}{11},\frac 32, Q^+(\varrho/1)}+\frac 1{\varrho}\| q^1\|_{\frac {12}{11},\frac 32, Q^+(\varrho/2)}\Big]\leq$$
$$\leq \Big[\frac 1{\varrho}\| \nabla v\|_{\frac {12}{11},\frac 32, Q^+(\varrho/2)}+\frac 1{\varrho}\|\nabla q\|_{\frac {12}{11},\frac 32, Q^+(\varrho/2)}+$$$$+\frac 1{\varrho}\| \nabla v^1\|_{\frac {12}{11},\frac 32, Q^+(\varrho/2)}+\frac 1{\varrho}\| q^1\|_{\frac {12}{11},\frac 32, Q^+(\varrho/2)}\Big].$$
So, by (\ref{weakerterm}), by (\ref{highder})
with $R=\varrho$ and $\tau=\frac 12$, and by (\ref{v1est}), one can find the following bound
\begin{equation}\label{q2est}
\|\nabla q^2\|_{9,\frac 32, Q^+(\varrho/4)}	\leq
\frac c{\varrho^\frac 43}\Big[E^\frac 12(v,\varrho)+D_2^\frac 23(q,\varrho)+$$$$+(A^\frac 38(\varrho)E^\frac 34(v,\varrho)C^\frac 14(v,\varrho))^\frac 23\Big].
\end{equation}
Now, assuming $0<r<\varrho/4$,  we can derive from (\ref{v1est}) and from (\ref{q2est}) the estimate
$$D_2(r)\leq \frac c{r^\frac {13}8}\Big[\int\limits^0_{-r^2}\Big(\int\limits_{B^+(r)}|\nabla q^1|^\frac {12}{11}dx\Big)^\frac {11}8dt+
\int\limits^0_{-r^2}\Big(\int\limits_{B^+(r)}|\nabla q^2|^\frac {12}{11}dx\Big)^\frac {11}8dt\Big]\leq 
$$$$
\leq \frac c{r^\frac {13}8}\int\limits^0_{-r^2}\Big(\int\limits_{B^+(r)}|\nabla q^1|^\frac {12}{11}dx\Big)^\frac {11}8dt+cr^2
\int\limits^0_{-r^2}\Big(\int\limits_{B^+(r)}|\nabla q^1|^9dx\Big)^\frac 16dt\leq$$
$$\leq c\Big(\frac \varrho r\Big)^\frac {13}8A^\frac 38(v,\varrho)E^\frac 34(v,\varrho)C^\frac 14(v,\varrho)+c\Big(\frac r\varrho\Big)^2\Big[E^\frac 34(v,\varrho)+D_2(q,\varrho)+$$$$+A^\frac 38(v,\varrho)E^\frac 34(v,\varrho)C^\frac 14(v,\varrho)\Big]
$$
and thus
$$D_2(q,r)\leq c\Big(\frac r\varrho\Big)^2\Big[E^\frac 34(v,\varrho)+D_2(q,\varrho)\Big]+$$$$+c\Big(\frac \varrho r\Big)^\frac {13}8A^\frac 38(v,\varrho)E^\frac 34(v,\varrho)C^14(v,\varrho)
$$
for $0<r<\varrho/4$. The latter implies estimate
(\ref{pressure}).

\section{Proof of Proposition \ref{boundednesstheorem}}
 \setcounter{equation}{0}

\begin{proof} We let $g=g(v)$ and $G=G_1(v,q)$.

Let us 
assume that $g<\infty$. Our aim is to show that $G<\infty$. There are three cases:

\textsc{Case 1}.
Suppose that
\begin{equation}\label{case1}
C_0:=\sup\limits_{0<R<1}C(v,R)<\infty.	
\end{equation}
Then, from (\ref{energy}), one can deduce  that
$$A(v, R/2)+E(v, R/2)\le c_1(1+D_0^\frac 23(q,R)).$$	
Here and in what follows in this case, $c_1$ is a generic constant depending on $C_0$ only. 

Now, let us use (\ref{embedding}), (\ref{pressure}) with $\varrho= R/2$, and the above estimate. As a result, we find
$$D_2(q,r)\leq c\Big(\frac rR\Big)^2D_2(q, R/2)
+c_1\Big(\frac Rr\Big)^\frac {13}8[E^\frac 34(v, R/2)+1 +D_2^\frac 34(q,R)]\leq 
$$$$\leq c\Big(\frac rR\Big)^2D_2(q,R)+c_1\Big(\frac Rr\Big)^\frac {13}8[1+D_2(q,R)^\frac 23]$$
for all $0<r< R/2$. So, by Young's inequality,
\begin{equation}\label{pressure1}
D_2(q,r)\leq c\Big(\frac rR\Big)^2D_2(q,R)+c_1\Big(\frac Rr\Big)^\frac {71}8\end{equation}
for all $0<r< R/2$. 
If $ R/2\leq r\leq R$,  then  
$$D_2(q,r)\leq \frac 1{( R/2)^\frac {13}8}\int\limits^0_{-R^2}\Big(\int\limits_{B^+(R)}|\nabla q|^\frac {12}{11}dx\Big)^\frac {11}8dt\leq 2^\frac {13}8D_2(q,R)\Big(\frac {2r}{R}\Big)^2.$$
So, estimate (\ref{pressure1})
 holds for all $0<r<R<1$. 

Now, for $\mu$ and  $R$ in $]0,1[$, we let $r=\mu R$ in  (\ref{pressure1}) and find
$$D_2(q,\mu R)\leq c\mu^2D_2(q,R)+c_1\mu^{-\frac{71}8}.
$$
Picking $\mu$ up so small that $2c\mu\leq 1$, we show that
$$D_2(q,\mu R)\leq \mu D_2(q,R)+c_1
$$ for any $0<R<1$. One can iterate the last inequality and get the following:
$$D_2(q,\mu^{k+1}R)\leq \mu^{k+1}D_2(q,R)+c_1(1+\mu+...+\mu^k)$$
for all natural numbers $k$. The latter implies 
that 
\begin{equation}\label{1case1est}
D_2(q,r)\leq c_1\frac rRD_2(q,R)+c_1	
\end{equation} for all $0<r<R<1$. And we can deduce from (\ref{embedding}) and from  the above estimate that
$$\max\{\sup\limits_{0<R<1}D_0(q,R), \sup\limits_{0<R<1}D_2(q,\tau R)\}<\infty$$
for any $0<\tau<1$. Uniform boundedness of $A(R)$ and $E(R)$ follows from the energy estimate (\ref{energy}) and from the assumption (\ref{case1}).

\textsc{Case 2}. Assume now that 
\begin{equation}\label{case2}
	A_0:=\sup\limits_{0<R<1}A(v,R)<\infty.
	\end{equation}

Then, from (\ref{multiple}), it follows that
$$C(v,r)\leq cA_0^\frac 34E^\frac 34(v,r)$$
for any $0<r<1$ and thus
$$A(v,\tau \varrho)+E(v,\tau \varrho)\leq c_3(A_0,\tau)\Big[E^\frac 12(v,\varrho)+E^\frac 14(v,\varrho)D_0^\frac 23(q,\varrho)+E^\frac 34(v,\varrho)\Big].
$$ for any $0<\tau<1$ and $0<\varrho<1$.

Our next step is an estimate for  the pressure quantity:
$$D_2(q,r)\leq c\Big(\frac r\varrho\Big)^2\Big[D_2(q,\varrho)+E^\frac34(v,\varrho)\Big]+c_2\Big(\frac \varrho r\Big)^\frac {13}8E^\frac {15}{16}(v,\varrho)\leq$$$$\leq  c\Big(\frac r\varrho\Big)^2D_2(q,\varrho)+c_2\Big(\frac \varrho r\Big)^\frac {13}8(E^\frac {15}{16}(v,\varrho)+1)$$ for any $0<r<\varrho<1$.
Here,  a generic constant, depending on $A_0$ only, is denoted by $c_2$.

Letting $r=\tau R$ and $\mathcal E(r):=A(v,r)+D_2(q,r)$, one can deduce  from latter inequalities, see also  (\ref{embedding}), the following estimates:
$$\mathcal E(\tau \varrho)\leq c\tau^2D_2(q,\varrho)+c_2\Big(\frac 1 \tau\Big)^\frac {13}8(E^\frac {15}{16}(v,\varrho)+1)+$$$$+c_3(A_0,\tau)\Big[E^\frac 12(v,\varrho)+E^\frac 14(v,\varrho)D_2^\frac 23(\varrho)+E^\frac 34(v,\varrho)\Big]\leq $$
$$\leq c\tau^2D_2(q,\varrho)+c_2\Big(\frac 1 \tau\Big)^\frac {13}8(E^\frac {15}{16}(v,\varrho)+1)+$$$$+c_3(A_0,\tau)\Big(\frac1{\tau}\Big)^4E^\frac 34(v,\varrho)+c_3(A_0,\tau)\Big[E^\frac 12(v,\varrho)+E^\frac 34(v,\varrho)\Big]\leq $$
$$\leq c\tau^2\mathcal E(\varrho)
+c_3(A_0,\tau).
$$
The rest of the proof is similar to what has been done in Case 1, see derivation of (\ref{1case1est}).

\textsc{Case 3}. Assume now that 
\begin{equation}\label{case3}
	E_0:=\sup\limits_{0<R<1}E(v,R)<\infty.
	\end{equation}

Indeed,  
$$C(v,r)\leq cE_0^\frac 34A^\frac 34(v,r)$$
for all $0<r\leq 1$. As to the pressure, we can find
$$D_2(\tau\varrho)\leq c\tau^2D_2(\varrho)+c_4(E_0,\tau)A^\frac {9}{16}(\varrho)
$$ for any $0<\tau<1$ and for any $0<\varrho<1$.
In turn, the energy inequality gives:
$$A(v,\tau \varrho)\leq c_5(E_0,\tau)\Big[A^\frac 12(v,\varrho)+A^\frac 14(v,\varrho)D_0^\frac 23(q,\varrho)+A^\frac 34(v,\varrho)\Big]\leq
$$
$$\leq c_5(E_0,\tau)\Big[A^\frac 12(v,\varrho)+A^\frac 14(v,\varrho)D_2^\frac 23(q,\varrho)+A^\frac 34(v,\varrho)\Big]
$$ for any $0<\tau<1$ and for any $0<\varrho<1$.
Similar to Case 2, one can introduce the quantity $\mathcal E(r)=A(v,r)+D_2(q,r)$ and find the following inequality for it:
$$\mathcal E(\tau\varrho)\leq  c\tau^2D_2(q,\varrho)+c_4(E_0,\tau)A^\frac {9}{16}(v,\varrho)+$$
$$+c_5(E_0,\tau)\Big[A^\frac 12(v,\varrho)+A^\frac 14(v,\varrho)D_2^\frac 23(q,\varrho)+A^\frac 34(v,\varrho)\Big]\leq $$
$$\leq c\tau^2\mathcal E(\varrho)+c_5(E_0,\tau)$$
for any $0<\tau<1$ and for any $0<\varrho<1$. The rest of the proof is the same as in Case 2. 

\end{proof}

\section{Proof of Theorem \ref{local energy ancient solution}}
 \setcounter{equation}{0}

Assume that $v$ and $q$ is a suitable weak solution in $Q^+$ with Type I blow up at the origin so that
\begin{equation}\label{type1}
g=g(v)=\min\{ \sup\limits_{0<R<1}A(v,R),\sup\limits_{0<R<1}E(v,R)\sup\limits_{0<R<1}C(v,R)\}<\infty.	
\end{equation}
By Theorem \ref{boundednesstheorem},
\begin{equation}
	\label{bound1}
G_1=G_1(v,q):=\max\{\sup\limits_{0<R<1}A(v,R),\sup\limits_{0<R<1}E(v,R),$$$$\sup\limits_{0<R<1}C(v,R),	\sup\limits_{0<R<1}D_0(v,R)\}<\infty.
\end{equation}
We know, see Theorem 2.2 in \cite{S2016}, that there exists a positive number $\varepsilon_1=\varepsilon_1(G_1)$ such that
\begin{equation}\label{sing1}
\inf\limits_{0<R<1}C(v,R)\geq \varepsilon_1>0.	
\end{equation}
Otherwise, the origin $z=0$ is a regular point of $v$.

Let $R_k\to0$ and $a>0$ and let
$$u^{(k)}(y,s)=R_kv(x,t),\qquad p^{(k)}(y,s)=R_k^2q(x,t),
$$ where $x=R_ky$, $t=R^2_ks$. Then, we have
$$A(v,aR_k)=A(u^{(k)},a)\leq G_1,\qquad E(v,aR_k)=E(u^{(k)},a)\leq G_1,$$$$ C(v,aR_k)=C(u^{(k)},a)\leq G_1, \qquad D_0(q,u^{(k)})=D_0(p^{(k)},a)\leq G_1.$$	
Thus, by (\ref{highder}),
$$\|\partial_tu^{(k)}\|_{\frac {12}{11},\frac 32,Q^+(a)}+\|\nabla^2u^{(k)}\|_{\frac {12}{11},\frac 32,Q^+(a)}+\|\nabla p^{(k)}\|_{\frac {12}{11},\frac 32,Q^+(a)}\leq c(a,G_1)$$
Moreover, the well known multiplicative inequality implies the following bound:
$$\sup\limits_k\int\limits_{Q^+}|u^{(k)}|^\frac {10}3dz\leq c(a,G_1).$$

Using known arguments,
one can select a subsequence (still denoted in the same way as the whole sequence) such that, for any $a>0$,
$$u^{(k)}\to u$$
in $L_3(Q^+(a))$,
$$\nabla u^{(k)}\rightharpoonup \nabla u$$
in $L_2(Q^+(a))$,
$$p^{(k)}\rightharpoonup p$$
in $L_\frac 32(Q^+(a))$.  The first two statements are well known and we shall comment on the last one only.

Without loss of generality, we may assume that
$$\nabla p^{(k)}\rightharpoonup w
$$ in $L_{\frac {12}{11}}(Q^+(a))$ for all positive $a$.

We let $p^{(k)}_1(x,t)=p^{(k)}(x,t)-[p^{(k)}]_{B^+(1)}(t)$. Then, there exists 
a subsequence $\{k^1_j\}_{j=1}^\infty$ such that 
$$p^{(k^1_j)}_1\rightharpoonup p_1$$ in $L_\frac 32(Q^+(1))$ as $j\to\infty$. Indeed, it follows from Poincar\'e-Sobolev inequality
$$\|p^{(k^1_j)}_1\|_{\frac 32, Q^+(1)}\leq c \|\nabla p^{(k^1_j)}\|_{\frac {12}{11},\frac 32,Q^+(1)}\leq c(1,G_1).$$
Moreover, one has $\nabla p_1=w$ in $Q^+(1)$. 

Our next step is to define $p^{(k^1_j)}_2(x,t)=p^{(k^1_j)}(x,t)-[p^{(k^1_j)}]_{B^+(2)}(t)$. For the same reason as above, there is a subsequence $\{k^2_j\}_{j=1}^\infty$ of the sequence 
$\{k^1_j\}_{j=1}^\infty$ such that 
$$p^{(k^2_j)}_2\rightharpoonup p_2$$ in $L_\frac 32(Q^+(2))$ as $j\to\infty$. Moreover, we claim that $\nabla p_2=w$ in $Q^+(2)$ and 
$$p_2(x,t)-p_1(x,t)=[p_2]_{B^+(1)}(t)-[p_1]_{B^+(1)}(t)=[p_2]_{B^+(1)}(t)$$
for $x\in B^+(1)$ and $-1<t<0$, i.e., in $Q^+(1)$.

After $s$ steps, we arrive at the following: there exists a subsequence $\{k^s_j\}_{j=1}^\infty$ of the sequence $\{k^{s-1}_j\}_{j=1}^\infty$ such that $p^{(k^s_j)}_s(x,t)
=p^{(k^s_j)}(x,t)-[p^{(k^s_j)}]_{B^+(s)}(t)$ in $Q^+(s)$ and 
$$p^{(k^s_j)}_s\rightharpoonup p_s$$ in $L_\frac 32(Q^+(s))$ as $j\to\infty$.
Moreover, $\nabla p_s=w$ in $Q^+(s)$ and 
$$p_s(x,t)=p_{s-1}(x,t)+[p_s]_{B^+(s-1)}(t)$$
in $Q^+(s-1)$. And so on.

The following function $p$  is going to be well defined: $p=p_1$ in $Q^+(1)$ and 
$$p(x,t)=p_{s+1}(x,t)-\sum\limits_{m=1}^s[p_{m+1}]_{B^+(m)}(t)\chi_{]-m^2,0[}(t)$$
in $Q^+(s+1)$, where $\chi_\omega(t)$  is the indicator function of the set $\omega \in \mathbb R$. Indeed, to this end, we need  to verify 
that 
$$p_{s+1}(x,t)-\sum\limits_{m=1}^{s}[p_{m+1}]_{B^+(m)}(t)\chi_{]-m^2,0[}(t)=$$
$$=p_{s}(x,t)-\sum\limits_{m=1}^{s-1}[p_{m+1}]_{B^+(m)}(t)\chi_{]-m^2,0[}(t)$$
in $Q^+(s)$. The latter is an easy exercise.  

Now,  let us fix $s$ and consider the sequence
$$p^{(k^s_j)}(x,t)=p^{(k^s_j)}_{s}(x,t)-\sum\limits_{m=1}^{s-1}[p^{(k^s_j)}_{m+1}]_{B^+(m)}(t)\chi_{]-m^2,0[}(t)$$
in $Q^+(s)$. Then, since the sequence $\{k^s_j\}_{j=1}^\infty$  is a subsequence of all sequences $\{k^{m+1}_j\}_{j=1}^\infty$ with $m\leq s-1$, one can easily check that
$$p^{(k^s_j)}\rightharpoonup p$$
in $L_\frac 32(Q^+(s))$. It remains to apply     the diagonal procedure of Cantor.

Having in hands the above convergences, we can conclude that the pair $u$ and $p$ is a local energy ancient solution in $Q^+_-$ and (\ref{scalequatities}) and (\ref{singularity}) hold.

The inverse statement is obvious.

\section{Proof of Theorem \ref{mbas_type1}}
\setcounter{equation}{0}

The proof is similar to the proof of Theorem \ref{local energy ancient solution}. We start with scaling $u^\lambda(y,s)=\lambda u(x,t)$ and $p^\lambda(y,s)=\lambda^2p(x,t)$ where $x=\lambda y$ and $t=\lambda^2s$ and $\lambda\to\infty$. We know
$$|u^\lambda(0,y_{3\lambda},0)|=\lambda|u(0,a,0)|=\lambda$$
and so that $y_{3\lambda}\to0$  as $\lambda\to\infty$.

For any $R>0$, by the invariance with respect to the scaling, we have
$$A(u^\lambda,R)=A(u,\lambda R)\leq G(u,p)=:G_0,\qquad E(u^\lambda,R)=E(u,\lambda R)\leq G_0,
$$
$$C(u^\lambda,R)=C(u,\lambda R)\leq G_0,\qquad D_0(p^\lambda,R)=D_0(p,\lambda R)\leq G_0.$$
Now, one can apply estimate (\ref{scalequatities}) and get the following:
$$\|\partial_tu^\lambda\|_{\frac {12}{11},\frac 32,Q^+(R)}+\|\nabla^2u^\lambda\|_{\frac {12}{11},\frac 32,Q^+(R)}+\|\nabla p^\lambda\|_{\frac {12}{11},\frac 32,Q^+(aR}\leq c(R,G_0).$$
Without loss of generality, we can deduce from the above estimates that, for any $R>0$,
$$u^{(k)}\to v$$
in $L_3(Q^+(R))$,
$$\nabla u^{(k)}\rightharpoonup \nabla v$$
in $L_2(Q^+(R))$,
$$p^{(k)}\rightharpoonup q$$
in $L_\frac 32(Q^+(R))$. Passing to the limit as $\lambda\to\infty$, we conclude that $v$ and $q$ are a local energy ancient solution in $Q^+_-$ for which $G(v,q)<\infty$.

Now, our goal is to prove that $z=0$ is a singular point of $v$. We argue ad absurdum. Assume that the origin is a regular point, i.e., there exist  numbers $R_0>0$ and $A_0>0$
such that 
$$|v(z)|\leq A_0$$
for all $z\in Q^+(R_0)$. Hence,
\begin{equation}\label{estim1}
C(v,R)=\frac 1{R^2}\int\limits_{Q^+(R)}|v|^3dz\leq cA_0^3R^3	
\end{equation}
for all $0<R\leq R_0$. Moreover,
\begin{equation}\label{pass}
	C(u^\lambda,R)\to C(v,R)
\end{equation}
as $\lambda\to\infty$. By weak convergence, 
$$D_0(q,R)\leq G_0$$
for all $R>0$.	
Now, we can calculate positive numbers $\varepsilon(G_0)$ and $c(G_0)$ of Theorem 2.2 in \cite{S2016}. Then, let us fix $0<R_1<R_0$, see (\ref{estim1}), so that $C(v,R_1)<\varepsilon(G_0)/2$. According to (\ref{pass}), one can find a number $\lambda_0>0$ such that 
$$G(u^\lambda,R_1)<\varepsilon(G_0)$$
for all $\lambda>\lambda_0$. By Theorem 2.2 of \cite{S2016}, 
$$\sup\limits_{z\in Q^+(R_1/2)}|u^\lambda(z)|<\frac {c(G_0)}{R_1}
$$ for all $\lambda>\lambda_0$.
It remains to select $\lambda_1>\lambda_0$ such that $y_{3\lambda}=a/\lambda <R_1/2$ and $\lambda_1>c(G_0)/R_1$. Then
$$|u^{\lambda_1}(0,y_{3\lambda_1},0)|=\lambda_1\leq \sup\limits_{z\in Q^+(R_1/2)}|u^{\lambda_1}(z)|<\frac {c(G_0)}{R_1}.$$
This is a contradiction.

\end{document}